\documentclass[12pt]{amsart}
\usepackage[margin=3cm, top=1in]{geometry}

\begin{document}
\title{Alternating permutations containing the pattern $321$ or $123$ exactly once}
\author{Joel Brewster Lewis (MIT)}
\maketitle

Inspired by a recent note of Zeilberger \cite{Z}, Alejandro Morales asked whether one can count alternating (i.e., up-down) permutations that contain the pattern $123$ or $321$ exactly once.  In this note we answer the question in the affirmative; in particular, we show that 
\[
a_{2m}(321) = \frac{4(m - 2)\cdot (2m + 3)!}{(m+1)! \cdot (m + 4)!}, \qquad a_{2m}(123) = \frac{10\cdot (2m)!}{(m-2)!\cdot (m+3)!}
\]
and
\[
a_{2m + 1}(321) = a_{2m + 1}(123) = \frac{3(3m +4)(m - 1) \cdot (2m + 2)!}{(m + 1)!\cdot (m + 4)!}
\]
where $a_{n}(p)$ is the number of alternating permutations of length $n$ containing the pattern $p$ exactly once.

Suppose an alternating permutation $w$ contains a single instance  of the pattern $321$, say $i < j < k$ and $w_i > w_j > w_k$.  
% If $j$ is odd then $w_{j - 1} > w_j > w_k$ is an instance of $321$ in $w$ and so in this case we must have $i = j - 1$, while if $j$ is even then $w_i > w_j > w_{j + 1}$ is an instance of $321$ in $w$ so in this case we must have $k = j + 1$.  
As in the Noonan/Burstein/Zeilberger \cite{N, B, Z} case, we conclude that $w$ contains a unique copy of $321$ if and only if the following conditions hold:
\begin{itemize}%[itemsep=-1ex, topsep=0ex]
\item other than $w_i$, the entries preceding $w_j$ in $w$ are less than $w_j$;
\item other than $w_k$, the entries following $w_j$ in $w$ are greater than $w_j$; and
\item the permutations $u = w_1w_2 \cdots w_{j - 1}w_k$ and $v = w_i w_{j + 1} w_{j + 2} \cdots w_n$ avoid $321$.
\end{itemize}
Under these conditions we have that 
\begin{enumerate}%[itemsep=-1ex, topsep=0ex]
\item $u$ is an alternating permutation of length $j$ not ending in its largest entry, and 
\item if $j$ is odd then $v$ is an alternating permutation not beginning with its smallest entry, while if $j$ is even then $v$ is a reverse-alternating (i.e., down-up) permutation not beginning with its smallest entry.
\end{enumerate}
Moreover, it's easy to check that an alternating permutation $w$ with a single instance of $321$ is uniquely determined by the choice of $321$-avoiding permutations order-isomorphic to $u$ and $v$ satisfying conditions 1 and 2.

Following the lead of Deutsch, Reifegerste and Stanley \cite[Problems $h^7$ and $i^7$]{catadd}, it is a straightforward exercise to show the enumerations shown in Table~\ref{table}.
\begin{table}
\begin{center}
\begin{tabular}{|c|c||c|c|c|c|}
\hline
type & length      & total number 									& ending in largest	 & beginning with smallest  & valid for \\ \hline
UD   & $2\ell$     & $C_{\ell + 1}$ & $C_\ell$ & $C_\ell$ & $\ell \geq 2$\\ \hline
UD   & $2\ell + 1$ & $C_{\ell + 1}$ & 0  										 & $C_\ell$ & $\ell \geq 1$ \\ \hline
DU   & $2\ell$     & $C_{\ell}$	& 0 										 & 0 & $\ell \geq 0$ \\ \hline
DU   & $2\ell + 1$ & $C_{\ell + 1}$ & $C_\ell$ & 0  & $\ell \geq 1$ \\ \hline
\end{tabular}
\end{center}
\caption{Enumerations of $321$-avoiding alternating permutations, as well as those ending in their largest entry or beginning with their smallest entry.}
\label{table}
\end{table}

It follows immediately that when $n = 2m$ ($m \geq 2$), the number of alternating permutations of length $n$ containing $321$ exactly once is
\[
\sum_{j = 1}^{m - 2} C_{j + 1}(C_{m - j + 1} - C_{m - j})  + \sum_{j = 2}^{m - 1} (C_{j + 1} - C_j)C_{m - j + 1} = \frac{4(m - 2)(2m + 3)!}{(m+1)!(m + 4)!},
\]
while when $n = 2m + 1$ ($m \geq 1$),  the number of alternating permutations of length $n$ containing $321$ exactly once is
\[
\sum_{j = 1}^{m - 1}C_{j + 1}(C_{m - j + 1} - C_{m - j})  + \sum_{j = 2}^m (C_{j + 1} - C_j)C_{m - j + 1} = \frac{3(3m +4)(m - 1) (2m + 2)!}{(m + 1)!(m + 4)!}.
\]
The number of alternating permutations of length $2m + 1$ containing $123$ exactly once is the same as the number of alternating permutations of length $2m + 1$ containing $321$ exactly once, since reversing any permutation of either type gives a permutation of the other.  And a similar analysis shows that for $m \geq 2$, the number of alternating permutations of length $2m$ containing $123$ exactly once is
\[
\frac{10(2m)!}{(m-2)!(m+3)!},
\]
and of course this is also the number of reverse-alternating permutations of length $2m$ containing $321$ exactly once.


\begin{thebibliography}{99}

\bibitem{B} Alex Burstein, A short proof for the number of permutations containing the pattern 321 exactly once, Elec. J. Comb. 18(2), \#P21, 2011.

\bibitem{N} John Noonan, The number of permutations containing exactly one increasing subsequence of length three, Discrete Math. 152, 307-313, 1996.

\bibitem{catadd} Richard P.\ Stanley, Catalan addendum to \emph{Enumerative Combinatorics, Vol.~2}, \texttt{http://www-math.mit.edu/$\sim$rstan/ec/catadd.pdf}, version of October 22, 2011.

\bibitem{Z} Doron Zeilberger, Alexander Burstein's lovely combinatorial proof of John Noonan's beautiful formula that the number of $n$-permutations that contain the pattern $321$ exactly once equals $(3/n)(2n)!/((n-3)!(n+3)!)$, \texttt{arXiv:1110.4379v2}, 2011.

\end{thebibliography}
\end{document}